# On pairs of complementary GJ pivoting transforming skew-symmetric matrices

Samuel Awoniyi
Department of Industrial and Manufacturing Engineering
FAMU-FSU College of Engineering
2525 Pottsdamer Street
Tallahassee, FL 32310
E-mail: awoniyi@eng.famu.fsu.edu
ORCID# 0000 0001 7102 6257

# Abstract

This article describes certain ratios that attend pairs of complementary Gauss-Jordan pivoting transforming skew-symmetric matrices. Our interest in those ratios was motivated by a need to prove a crucial Claim stated in a recently proposed strongly polynomial-time algorithm for the general linear programming problem. That Claim is proved in this article and, as a consequence of this proof, a compact implementation of the strongly polynomial-time algorithm is outlined.

# 1. Introduction

This article describes certain ratios that emerge when pairs of Gauss-Jordan (GJ) pivoting transform skew-symmetric matrices. Our motivation for this article is a need to give a proof of Claim 4 in the strongly polynomial-time algorithm described in [1] for the general linear programming (LP) problem.

Claim 4 in [1] is about certain ratios being equal to one another when certain pairs of GJ pivotings are implemented as iterations of the proposed algorithm. The article [1] is available as arXiv:2503.12041 [pdf] at www.arxiv.org.

The remainder of this article is organized as follows. Section 2 gives a numerical example that illustrates our practical motivation for this article. Section 3 states and proves a new basic lemma regarding how pairs of (complementary) GJ pivoting transform skew-symmetric matrices. Section 4 applies the lemma of Section 3 to give a proof of Claim 4 of the article described in [1].

# 2. A motivational example

In this Section, we will show how the strongly polynomial-time algorithm proposed in [1] solves a simple LP problem, and thereby indicate the origin of our interest in how pairs of GJ pivotings transform skew-symmetric matrices.

In this article, the general LP problem is assumed to be given in Neumann symmetric form, (P) below:

.

$$\left\{ \begin{array}{ll} \text{maximize} & f^T x \\ \text{subject to:} & Ax \leq b, \\ & x \geq 0 \end{array} \right\} \cdots\cdots (P)$$

.

where $f$ is $n$-vector, $A$ is $k$-by-$n$ matrix, $b$ is $k$-vector, and $x$ is $n$-vector of the problem's variables. From basic LP duality theory, solving (P) is equivalent to computing a $2(k+n)$-vector $z$ that solves the constrained system of linear equations (Eq) stated below:

.

$$\left\{\begin{array}{l} Mz = q, \\ z_j z_{(k+n+j)} = 0, \ \text{for}\, j = 1,\ldots,k+n. \\ z \geq 0 \end{array}\right\} \cdots\cdots(\text{Eq})$$

.

where

.

$$M = \begin{pmatrix} O & A & I_{(k)} & O \\ -A^T & O & O & I_{(n)} \\ -b^T & f^T & o^T & o^T \end{pmatrix} \text{ and } q = \begin{pmatrix} b \\ -f \\ o \end{pmatrix}$$

To illustrate that problem definition and the iterations of the algorithm proposed in [1] for solving (Eq), we consider the following LP problem

$$\begin{array}{rl} \text{maximize:} & 2x_1 + x_2 \\ \text{such that:} & x_1 + x_2 \leq 5 \\ & x_1 \leq 2 \\ & x_1 \geq 0,\ x_2 \geq 0 \end{array}$$

Here,

$$A = \begin{pmatrix} 1 & 1 \\ 1 & 0 \end{pmatrix},\ b = \begin{pmatrix} 5 \\ 2 \end{pmatrix} \ \&\, f = \begin{pmatrix} 2 \\ 1 \end{pmatrix}.$$

so that

$$M = \begin{array}{|c|c|c|c|c|c|c|c|} \hline 0 & 0 & 1 & 1 & 1 & 0 & 0 & 0 \\ \hline 0 & 0 & 1 & 0 & 0 & 1 & 0 & 0 \\ \hline -1 & -1 & 0 & 0 & 0 & 0 & 1 & 0 \\ \hline -1 & 0 & 0 & 0 & 0 & 0 & 0 & 1 \\ \hline -5 & -2 & 2 & 1 & 0 & 0 & 0 & 0 \\ \hline \end{array} \quad \&\ q = \begin{array}{|c|} \hline 5 \\ \hline 2 \\ \hline -2 \\ \hline -1 \\ \hline 0 \\ \hline \end{array}$$

To utilize the same notation used in [1], let $[M\ q]$ denote the augmented matrix formed from $M$ and $q$, and let "$[M\ q]$ instance" denote the result of transforming $[M\ q]$ with some elementary row operations such as the ones that comprise GJ pivoting.

In the algorithm proposed in [1], (Eq) is solved by first *adding the last row to every other row of* $[M\ q]$ and, thereafter, by *applying pairs of complementary GJ pivoting* until $q \geq 0$ with $q_{k+n+1} = 0$ in an $[M\ q]$ instance.

Applying the algorithm described in [1] to the particular example stated above, we see that the first $[M\ q]$ instance, at "initialization" is P0:

P0 =

| -5 | -2 | 3 | 2 | 1 | 0 | 0 | 0 | 5 |
|---|---|---|---|---|---|---|---|---|
| -5 | -2 | 3 | 1 | 0 | 1 | 0 | 0 | 2 |
| -6 | -3 | 2 | 1 | 0 | 0 | 1 | 0 | -2 |
| -6 | -2 | 2 | 1 | 0 | 0 | 0 | 1 | -1 |
| -5 | -2 | 2 | 1 | 0 | 0 | 0 | 0 | 0 |

.

After two GJ pivotings, namely, the first GJ pivoting at (4,4), followed by the second GJ pivoting at (1,1), the resulting $[M\ q]$ instance is P1:

P1 =

| 1 | 0.29 | -0.14 | 0 | 0.14 | 0 | 0 | -0.29 | 1 |
|---|---|---|---|---|---|---|---|---|
| 0 | -0.29 | 1.14 | 0 | -0.14 | 1 | 0 | -0.71 | 2 |
| 0 | -1 | 0 | 0 | 0 | 0 | 1 | -1 | -1 |
| 0 | -0.29 | 1.14 | 1 | 0.86 | 0 | 0 | -0.71 | 5 |
| 0 | -0.29 | 0.14 | 0 | -0.14 | 0 | 0 | -0.71 | 0 |

.

Next, after two more GJ pivotings, namely, the first GJ pivoting at (3,3), followed by the second GJ pivoting at (2,2), the resulting $[M\ q]$ instance is P2:

P2 =

| 1 | 0 | 0 | 0 | 0.1 | 0.1 | 0.2 | -0.7 | 1 |
|---|---|---|---|---|---|---|---|---|
| 0 | 1 | 0 | 0 | 0.1 | 0.1 | -0.8 | 1.3 | 1 |
| 0 | 0 | 1 | 0 | -0.1 | 0.9 | -0.2 | -0.3 | 2 |
| 0 | 0 | 0 | 1 | 1 | -1 | 0 | 0 | 3 |
| 0 | 0 | 0 | 0 | -1 | -1 | -2 | -3 | 0 |

.

The algorithm terminates with the $[M\ q]$ instance labelled P2, because a solution of corresponding (Eq) is reached there. Let us next note certain ratios in P0, P1 & P2.

.

**P0, along with a last-column-last-row correspondence enumeration, is

| -5 | -2 | 3 | 2 | 1 | 0 | 0 | 0 | 5 | ←1 |
|---|---|---|---|---|---|---|---|---|---|
| -5 | -2 | 3 | 1 | 0 | 1 | 0 | 0 | 2 | ←2 |
| -6 | -3 | 2 | 1 | 0 | 0 | 1 | 0 | -2 | ←3 |
| -6 | -2 | 2 | 1 | 0 | 0 | 0 | 1 | -1 | ←4 |
| -5 | -2 | 2 | 1 | 0 | 0 | 0 | 0 | 0 | |
| ↑ | ↑ | ↑ | ↑ | | | | | | |
| 1 | 2 | 3 | 4 | | | | | | |

In P0, the $m_{5,i}/q_i$ ratio is $-1$, for $i = 1,2,3,4$. We will demonstrate later in this article that the equality among those four ratios is not accidental.

.

**P1, along with a last-column-last-row correspondence enumeration, is

| 1 | 0.29 | -0.14 | 0 | 0.14 | 0 | 0 | -0.29 | 1 | ←1 |
|---|---|---|---|---|---|---|---|---|---|
| 0 | -0.29 | 1.14 | 0 | -0.14 | 1 | 0 | -0.71 | 2 | ←2 |
| 0 | -1 | 0 | 0 | 0 | 0 | 1 | -1 | -1 | ←3 |
| 0 | -0.29 | 1.14 | 1 | 0.86 | 0 | 0 | -0.71 | 5 | ←4 |
| 0 | -0.29 | 0.14 | 0 | -0.14 | 0 | 0 | -0.71 | 0 | |
| | ↑ | ↑ | | ↑ | | | ↑ | | |
| | 2 | 3 | | 1 | | | 4 | | |

In P1, the $m_{5,i}/q_i$ (or else $m_{5,4+i}/q_i$) ratio here is -0.14, for $i = 1,2,3,4$. Again, we will demonstrate later in this article that the equality among those four ratios is not accidental.

.

**P2, along with a last-column-last-row correspondence enumeration, is

| 1 | 0 | 0 | 0 | 0.1 | 0.1 | 0.2 | -0.7 | 1 | ←1 |
|---|---|---|---|---|---|---|---|---|---|
| 0 | 1 | 0 | 0 | 0.1 | 0.1 | -0.8 | 1.3 | 1 | ←2 |
| 0 | 0 | 1 | 0 | -0.1 | 0.9 | -0.2 | -0.3 | 2 | ←3 |
| 0 | 0 | 0 | 1 | 1 | -1 | 0 | 0 | 3 | ←4 |
| 0 | 0 | 0 | 0 | -1 | -1 | -2 | -3 | 0 | |
| | | | | ↑ | ↑ | ↑ | ↑ | | |
| | | | | 1 | 2 | 3 | 4 | | |

In P2, the $m_{5,4+i}/q_i$ ratio here is -1, for $i = 1,2,3,4$. Thus, a certain *ratio-equality property* holds in each one of P0, P1 & P2. In Lemma 1 below, we will generalize that ratio-equality property.

To set the stage for Lemma 1, we note here one more thing about P0, P1 & P2 - that each one of P0, P1 & P2 contains exactly four unit-vector columns, and five non-unit-vector columns. The submatrices formed by non-unit-vector columns inside P0, P1 & P2 are 5-by-5 matrices, and are respectively

| -5 | -2 | 3 | 2 | 5 |
|---|---|---|---|---|
| -5 | -2 | 3 | 1 | 2 |
| -6 | -3 | 2 | 1 | -2 |
| -6 | -2 | 2 | 1 | -1 |
| -5 | -2 | 2 | 1 | 0 |

,

| 0.14 | 0.29 | -0.14 | -0.29 | 1 |
|---|---|---|---|---|
| -0.14 | -0.29 | 1.14 | -0.71 | 2 |
| 0 | -1 | 0 | -1 | -1 |
| 0.86 | -0.29 | 1.14 | -0.71 | 5 |
| -0.14 | -0.29 | 0.14 | -0.71 | 0 |

. &

.

| 0.1 | 0.1 | 0.2 | -0.7 | 1 |
|---|---|---|---|---|
| 0.1 | 0.1 | -0.8 | 1.3 | 1 |
| -0.1 | 0.9 | -0.2 | -0.3 | 2 |
| 1 | -1 | 0 | 0 | 3 |
| -1 | -1 | -2 | -3 | 0 |

For general LP problems as stated above, each non-unit-vector submatrix is a (k+n+1)-by-(k+n+1) matrix. Lemma 1 is about the ratio-equality property inside such (k+n+1)-by-(k+n+1) matrices.

# 3. A lemma on GJ pairs transforming skew-symmetric matrices

In this Section, we will begin by defining and illustrating the concept of "latently-skew-symmetric" matrix and the concept of "complementary Gauss-Jordan-plus (GJ+) pivoting in column *j*" of a square matrix. Thereafter, we will state a basic lemma (Lemma 1) on "pairs of GJ+ pivoting transforming skew-symmetric matrices". After illustrating Lemma 1, we will give a proof by induction.

.

DEFINITION 1: *An s-by-s matrix, say S, is called latently-skew-symmetric (to be abbreviated as L-skew-sym) if, for* $i = 1,\ldots,s$, *there exists a number* $\beta_i$ *such that adding* $\beta_i$ *multiple of row s (last row) to row i transforms S into a skew-symmetric matrix. Clearly,* $S_{s,s}$ *must be 0. As a consequence of this definition, the ratio* $S_{s,i}/S_{i,s}$ *is the same number for* $i = 1,\ldots,s{-}1$ *having* $S_{i,s} \neq 0$. *We will sometimes refer to that common ratio as the "last-row/last-column" ratio.*

.

Illustration of definition 1: *S* below is L-skew-sym, with last-row/last-column ratio equal to -2.

*S* =

| 2 | -8 | 1 |
|---|---|---|
| 4 | -6 | -3 |
| -2 | 6 | 0 |

.

To see that *S* is L-skew-sym, set $\beta_1 = 1,\ \beta_2 = 1,\ \beta_3 = -0.5$ (in the definition above), and then have the skew symmetric matrix

.

| 0 | -2 | 1 |
|---|---|---|
| 2 | 0 | -3 |
| -1 | 3 | 0 |

.

DEFINITION 2: *Given an s-by-s matrix S, along with a column of S, say column j, a "complementary Gauss-Jordan-plus (GJ+) pivoting in column j" is defined as the transformation of S that takes S and j as input, and then utilizes operations (i), (ii), (iii) below to return the matrix denoted by* $GJ_j^+(S)$ *as its output:*

*(i) "attach" the j-th unit vector to S and thereby have an s-by-(s+1) augmented matrix, say S*1, *whose (s+1)-th column is the j-th unit vector;*

*(ii) perform GJ pivoting at position (j,j) in the matrix S*1 *and label the resultant matrix as S*2;

*(iii) swap column j and column s+1 of S*2, *and thereafter drop the resultant column s+1 (which is the j-th unit vector); then label the resultant s-by-s matrix as* $GJ_j^+(S)$.

.

Illustration of definition 2: Let us again use matrix

$S =$

| 2 | -8 | 1 |
|---|---|---|
| 4 | -6 | -3 |
| -2 | 6 | 0 |

with $j = 2$. Then we have augmented matrix

$S1 =$

| 2 | -8 | 1 | 0 |
|---|---|---|---|
| 4 | -6 | -3 | 1 |
| -2 | 6 | 0 | 0 |

;

.

After applying GJ with pivot at (2,2), we have

$S2 =$

| -3.33 | 0 | 5 | -1.33 |
|---|---|---|---|
| -0.67 | 1 | 0.5 | -0.17 |
| 2.00 | 0 | -3 | 1 |

;

.

After swapping column 2 and column 4, and thereafter dropping off new column 4 (the 2-unit vector), we have

$GJ_2^+(S) =$

| -3.33 | -1.33 | 5 |
|---|---|---|
| -0.67 | -0.17 | 0.5 |
| 2.00 | 1 | -3 |

.

LEMMA 1: *Suppose S is an s-by-s L-skew-sym matrix, with* $s \geq 5$. *If necessary, add row s to another row of S, in order to ensure that the resultant matrix is an s-by-s L-skew-sym matrix that has two columns, say column* $j_1$ *and column* $j_2$, *wherein two consecutive complementary GJ+ pivoting are feasible.*

*Applying complementary GJ+ pivoting in column* $j_1$, *followed by applying GJ+ pivoting in column* $j_2$, *transforms S into an L-skew-sym matrix.*

.

We will give a proof by induction, after an illustration is presented to foster some intuition.

.

`Two short pre-proof examples to illustrate Lemma 1`

.

*First illustration: a numerical example* – We utilize the L-skew-sym matrix

$S =$

| 2 | -8 | 1 |
|---|---|---|
| 4 | -6 | -3 |
| -2 | 6 | 0 |

Let $j_1 = 2$ and $j_2 = 1$. Recall from "Illustration of Definition 2" above that

$GJ_2^+(S) =$

| -3.33 | -1.33 | 5 |
|---|---|---|
| -0.67 | -0.17 | 0.5 |
| 2.00 | 1 | -3 |

so that

$GJ_1^\dagger(GJ_2^\dagger(S))$ =

| -0.3 | 0.4 | -1.5 |
|---|---|---|
| -0.2 | 0.1 | -0.5 |
| 0.6 | 0.2 | 0 |

.

.

To see that $GJ_1^\dagger(GJ_2^\dagger(S))$ is L-skew-sym, set $\beta_1 = 0.5$, $\beta_2 = -0.5$, $\beta_3 = 1.5$ (in the definition of L-skew-sym above), and then have the skew-symmetric matrix

| 0 | 0.5 | -1.5 |
|---|---|---|
| -0.5 | 0 | -0.5 |
| 1.5 | 0.5 | 0 |

.

*Second illustration: a 3-by-3 general example* – We utilize the skew-symmetric matrix

| 0 | $S_{12}$ | $q_1$ |
|---|---|---|
| $-S_{12}$ | 0 | $q_2$ |
| $-q_1$ | $-q_2$ | 0 |

After adding row 3 to rows 1 &2, that becomes the L-skew-sym matrix

S =

| $-q_1$ | $S_{12}-q_2$ | $q_1$ |
|---|---|---|
| $-S_{12}-q_1$ | $-q_2$ | $q_2$ |
| $-q_1$ | $-q_2$ | 0 |

Then

$GJ_1^\dagger(S)$ =

| $-1/q_1$ | $-(S_{12}-q_2)/q_1$ | -1 |
|---|---|---|
| $-(S_{12}+q_1)/q_1$ | $S_{12}(q_2-q_1-S_{12})/q_1$ | $q_2-q_1-S_{12}$ |
| -1 | $-S_{12}$ | $-q_1$ |

and

.

$GJ_2^\dagger(GJ_1^\dagger(S))$ =

| $q_2/(S_{12}(q_2-q_1-S_{12}))$ | $(S_{12}-q_2)/(S_{12}(q_2-q_1-S_{12}))$ | $-q_2/S_{12}$ |
|---|---|---|
| $-(S_{12}+q_1)/(S_{12}(q_2-q_1-S_{12}))$ | $q_1/(S_{12}(q_2-q_1-S_{12}))$ | $q_1/S_{12}$ |
| $-q_2/(q_2-q_1-S_{12})$ | $q_1/(q_2-q_1-S_{12})$ | 0 |

.

To see that $GJ_2^\dagger(GJ_1^\dagger(S))$ is L-skew-sym, set $\beta_1 = 1/S_{12}$; $\beta_2 = -1/S_{12}$, $\beta_3 = (q_1-q_2)/S_{12}$, and then have the skew-symmetric matrix.

| 0 | $-1/S_{12}$ | $-q_2/S_{12}$ |
|---|---|---|
| $1/S_{12}$ | 0 | $q_1/S_{12}$ |
| $q_2/S_{12}$ | $-q_1/S_{12}$ | 0 |

.

## Proof of Lemma 1

.

This is a proof by induction on the dimension $s$ of $S$. The "Induction base case" is the case wherein $s = 4$. The "Induction step" will assume that the result is true for dimension $s = r$, and then demonstrate that the result is true for dimension $s = r + 1$ as well.

.

*Induction base case – S an L-skew-sym matrix with dimension s=4.*

.

For the "Induction base case", we will use $s = 4$, and we start with the general skew-symmetric matrix

| 0 | $S_{12}$ | $S_{13}$ | $q_1$ |
|---|---|---|---|
| $-S_{12}$ | 0 | $S_{23}$ | $q_2$ |
| $-S_{13}$ | $-S_{23}$ | 0 | $q_3$ |
| $-q_1$ | $-q_2$ | $-q_3$ | 0 |

We obtain S from that matrix by adding row 4 to row 3. We then do the first GJ+ with pivot at position (3,3) and the second GJ+ with pivot at position (2,2), and have

.

$GJ_2^+(GJ_3^+(S)) =$

| $\phi S_{12}/\pi_2 S_{23}$ | $(q_3 S_{12}-\pi_2 S_{13})/\pi_2 S_{23}$ | $S_{12}/\pi_2$ | $-\phi/S_{23}$ |
|---|---|---|---|
| $(q_3 S_{12}+\pi_1 S_{23})/\pi_2 S_{23}$ | $-q_3/\pi_2 S_{23}$ | $-1/\pi_2$ | $-q_3/S_{23}$ |
| $-S_{12}/S_{23}$ | $1/S_{23}$ | 0 | $(\pi_2-S_{23})/S_{23}$ |
| $\phi/\pi_2$ | $q_3/\pi_2$ | $(S_{23}-\pi_2)/\pi_2$ | 0 |

.

where $\pi_1 = S_{13}+q_1$; $\pi_2 = S_{23}+q_2$; and $\phi = S_{13}\pi_2 - S_{23}\pi_1 - S_{12}q_3$. To see that this $GJ_2^+(GJ_3^+(S))$ is L-skew-sym, set $\beta_1 = -S_{12}/S_{23}$, $\beta_2 = 1/S_{23}$, $\beta_3 = 0$ & $\beta_4 = \pi_2/S_{23}$, and have the skew-symmetric matrix

| 0 | $-S_{13}/S_{23}$ | $S_{12}/S_{23}$ | $-\phi/S_{23}$ |
|---|---|---|---|
| $S_{13}/S_{23}$ | 0 | $-1/S_{23}$ | $-q_3/S_{23}$ |
| $-S_{12}/S_{23}$ | $1/S_{23}$ | 0 | $(\pi_2-S_{23})/S_{23}$ |
| $\phi/S_{23}$ | $q_3/S_{23}$ | $-(\pi_2-S_{23})/S_{23}$ | 0 |

.

*Induction step - assume that the result is true for dimension $s = r$ and then demonstrate that it is true for dimension $s = r + 1$ as well.*

.

We assume the result is true for $s = r$, and then consider $s = r + 1$. We will begin by making several background notes, Notes (i), (ii) & (iii) below, that will be needed later in this part of the proof.

.

`Note (i):` If column 1 and column 2 are swapped in a skew-symmetric matrix, say T, of dimension at least 3, together with row 1 and row 2 being swapped at the same time, then the resultant matrix, denoted by $T^{(1\leftrightarrow 2)}$, is skew-symmetric as well. The following matrices illustrate that statement.

$T =$

| 0 | -3 | -2 |
|---|---|---|
| 3 | 0 | -4 |
| 2 | 4 | 0 |

& $T^{(1\leftrightarrow 2)} =$

| 0 | 3 | -4 |
|---|---|---|
| -3 | 0 | -2 |
| 4 | 2 | 0 |

.

`Note (ii):` For any skew-symmetric matrix, say T, of dimension t, the (t-1)-by-(t-1) submatrix defined by "columns 2, . . ., t, together with rows 2, . . ., t", which we denote as $T_{\{2\to t\}}$, is skew-symmetric as well. The following table illustrates that statement.

$T =$

| 0 | $T_{1,2} \dots T_{1,t}$ |
|---|---|
| $T_{2,1}$ ⋮ $T_{t,1}$ | $T_{\{2\to t\}}$ |

By virtue of the induction hypothesis, applying a given pair of GJ+ pivotings to $T_{\{2\to t\}}$ results in a (t-1)-by-(t-1) L-skew-sym matrix; let $G(T_{\{2\to t\}})$ denote the skew-symmetric matrix obtained from that L-skew-sym matrix (through adding a suitable multiple of the last row to every row).

.

`Note (iii):` For any t-by-t skew-symmetric matrix, say T, let $\Gamma(G(T_{\{2\to t\}})$ denote the t-by-t matrix obtained from $G(T_{\{2\to t\}})$ by (a) appending an all-zero column vector as first column, and (b) appending an all-zero row vector as first row.

The following matrices are intended to illustrate Notes (i), (ii) & (iii). Let

T =

| 0 | 1 | 2 | 3 | 4 |
|---|---|---|---|---|
| -1 | 0 | 5 | 6 | 7 |
| -2 | -5 | 0 | 8 | 9 |
| -3 | -6 | -8 | 0 | 10 |
| -4 | -7 | -9 | -10 | 0 |

; $T_{\{2\to 5\}} =$

| 0 | 5 | 6 | 7 |
|---|---|---|---|
| -5 | 0 | 8 | 9 |
| -6 | -8 | 0 | 10 |
| -7 | -9 | -10 | 0 |

;

.

$G(T_{\{2\to 5\}}) =$

| 0 | -0.75 | 0.625 | 6.5 |
|---|---|---|---|
| 0.75 | 0 | -0.125 | -1.25 |
| -0.625 | 0.125 | 0 | 1.125 |
| -6.5 | 1.25 | -1.125 | 0 |

;

.

$\Gamma(G(T_{\{2\to5\}}) =$

| | | | | |
|---|---|---|---|---|
| 0 | 0 | 0 | 0 | 0 |
| 0 | 0 | -0.75 | 0.625 | 6.5 |
| 0 | 0.75 | 0 | -0.125 | -1.25 |
| 0 | -0.625 | 0.125 | 0 | 1.125 |
| 0 | -6.5 | 1.25 | -1.125 | 0 |

.

For a comparison, the skew-symmetric matrix corresponding to the L-skew-symm obtained from applying the same pair of GJ+ to T, denoted as $G(T)$, is:

G(T)=

| | | | | |
|---|---|---|---|---|
| 0 | 1.375 | -0.375 | 0.25 | 3.125 |
| -1.375 | 0 | -0.75 | 0.625 | 6.5 |
| 0.375 | 0.75 | 0 | -0.125 | -1.25 |
| -0.25 | -0.625 | 0.125 | 0 | 1.125 |
| -3.125 | -6.5 | 1.25 | -1.125 | 0 |

When one replaces the matrix $T_{\{2\to5\}}$ with the matrix $(T^{(1\leftrightarrow2)})_{\{2\to5\}}$, one obtains $\Gamma(G((T^{(1\leftrightarrow2)})_{\{2\to5\}})$. One can thereafter do "a swapping of columns 1&2, together with a swapping of rows 1&2" of $\Gamma(G((T^{(1\leftrightarrow2)})_{\{2\to5\}})$ for a comparison of non-zero entries of $(\Gamma(G((T^{(1\leftrightarrow2)})_{\{2\to5\}}))^{(1\leftrightarrow2)}$ and G(T).

Similarly, when one replaces the matrix $T_{\{2\to5\}}$ with the matrix $(T^{(1\leftrightarrow3)})_{\{2\to5\}}$, one can do "a swapping of columns 1&3, together with a swapping of row 1&3" of $\Gamma(G((T^{(1\leftrightarrow3)})_{\{2\to5\}}$ to enable a comparison of non-zero entries of $(\Gamma(G((T^{(1\leftrightarrow3)})_{\{2\to5\}}))^{(1\leftrightarrow3)}$ and G(T).

`Induction step continued:` We now return to demonstrating the "Induction step"; this will consist of generalizing the illustration example matrices above where t=5.

For this "Induction step", suppose we have an L-skew-sym matrix, say T, of dimension t. One can assume without loss of generality that T is skew-symmetric.

Utilizing concepts introduced in the background notes above, to complete this induction step, we only need to show that the entries in the three (t–1)-by-(t-1) skew-symmetric matrices, $T_{\{2\to t\}}$, $(T^{(1\leftrightarrow2)})_{\{2\to t\}}$ & $(T^{(1\leftrightarrow3)})_{\{2\to t\}}$ suitably include all the entries contained in the t-by-t skew-symmetric matrix T.

To that end, it is clear that: the entries in the matrix $T_{\{2\to t\}}$ include all entries in T, except the entries in the first row of T and the first column of T; the entries in the first row and the last column of matrix $(T^{(1\leftrightarrow2)})_{\{2\to t\}}$ include all the entries in the first row of T and the first column of T, except two entries (in positions (1,2) and (2,1)) which the matrix $(T^{(1\leftrightarrow3)})_{\{2\to t\}}$ includes.

With t=5, the following is an illustration of the statement about including the entries of T. Let

T=

| | | | | |
|---|---|---|---|---|
| 0 | 1 | 2 | 3 | 4 |
| -1 | 0 | 5 | 6 | 7 |
| -2 | -5 | 0 | 8 | 9 |
| -3 | -6 | -8 | 0 | 10 |
| -4 | -7 | -9 | -10 | 0 |

; $T_{\{2\to5\}} =$

| | | | |
|---|---|---|---|
| 0 | 5 | 6 | 7 |
| -5 | 0 | 8 | 9 |
| -6 | -8 | 0 | 10 |
| -7 | -9 | -10 | 0 |

;

.

$(T^{(1\leftrightarrow 2)})_{\{2\to 5\}} =$

| 0 | 2 | 3 | 4 |
|---|---|---|---|
| -2 | 0 | 8 | 9 |
| -3 | -8 | 0 | 10 |
| -4 | -9 | -10 | 0 |

; $(T^{(1\leftrightarrow 3)})_{\{2\to 5\}} =$

| 0 | -1 | 6 | 7 |
|---|---|---|---|
| 1 | 0 | 3 | 4 |
| -6 | -3 | 0 | 10 |
| -7 | -4 | -10 | 0 |

.

As demonstrated in the clarification notes stated above, applying any given pair of GJ+ pivoting (to T, $T_{\{2\to 5\}}$, $(T^{(1\leftrightarrow 2)})_{\{2\to 5\}}$ & $(T^{(1\leftrightarrow 3)})_{\{2\to 5\}}$) is straightforward, and that concludes this induction step. ■

# 4. A proof of Claim 4 of [1]

Let's begin by recalling Claim 4 of [1].

CLAIM 4: *In Section 4.3, under MinorP pivoting instance description (that is, the* $[M\ q]$ *instance having* $q_{k+n+1} = 0$*), the ratio* $m_{k+n+1,i}\,/q_i$ *(or, else,* $m_{k+n+1,k+n+i}/q_i$*, when* $m_{k+n+1,i} = 0$*) has the same value for all row indices* $i = 1,..,k+n$ *having* $q_i \neq 0$ *in that* $[M\ q]$ *instance.*

.

PROOF OF CLAIM 4

.

This proof consists of showing that Claim 4 is a particular application of Lemma 1 above. To that end, let us recall Lemma 1 here, for ease of referencing:

.

*"Lemma 1: Suppose S is an s-by-s L-skew-sym matrix, with* $s \geq 5$*. If necessary, add row s to another row of S, in order to ensure that the resultant matrix is an s-by-s L-skew-sym matrix that has two columns, say column* $j_1$ *and column* $j_2$*, wherein two consecutive complementary GJ+ pivoting are feasible."*

*Applying complementary GJ+ pivoting in column* $j_1$*, followed by applying GJ+ pivoting in column* $j_2$*, transforms S into an L-skew-sym matrix.*

.

Recall that, in Section 2 (A Motivational Example), an equivalence relation/connection is established between the (k+n+1)-by-(2k+2n+1) matrices $[M\ q]$ utilized by the algorithm proposed in [1] and the (k+n+1)-by-(k+n+1) matrices utilized by Lemma 1. The core of that connection is that each $[M\ q]$ instance (that has $q_{k+n+1} = 0$) contains exactly k+n columns that are unit vectors, and exactly k+n+1 other columns that form an L-skew-symm matrix.

The pair of GJ pivoting on $[M\ q]$ (in each iteration of the algorithm of [1]) translates into the pair of GJ+ described in Lemma 1. Accordingly, S being an L-skew-sym matrix in the conclusion of Lemma 1 directly translates into "*the ratio* $m_{k+n+1,i}\,/q_i$ *(or, else,* $m_{k+n+1,k+n+i}/q_i$*, when* $m_{k+n+1,i} = 0$*) has the same value for all row indices* $i = 1,..,k+n$ *having* $q_i \neq 0$ *in that* $[M\ q]$ *instance.*" in the conclusion of Claim 4 above.■

This proof of Claim 4 suggests a compact implementation of the algorithm described in [1] wherein each iteration will use a (k+n+1)-by-(k+n+1) matrix in place of a (k+n+1)-by-(2k+2n+1) matrix.

# 5. Some directions for further work

One obvious direction for further work is to do some numerical experiments with the compact implementation suggested by this article. This compact implementation could make the LP algorithm proposed in [1] well-suited to some classes of real-world mixed-integer LP problems.

Lemma 1 above is a new information for the literature on GJ pivoting transforming skew-symmetric matrices. This may be ordinarily useful information for mathematical modeling, as the matrix $A-A^T$ is skew-symmetric for any square numerical matrix A [5].

.